%% file: Main.tex
\documentclass[review]{elsarticle}

\usepackage{lineno,hyperref}

\usepackage{todonotes}
\usepackage{makecell}
\usepackage{graphicx}
\usepackage{adjustbox}
\usepackage{paralist}
\usepackage{cleveref}
\usepackage{amsmath}
\usepackage{booktabs}
\usepackage{tabularx}
\modulolinenumbers[5]

\newcolumntype{R}[1]{>{\tiny\raggedleft\arraybackslash}b{#1}}

\journal{Journal of \LaTeX\ Templates}
\input{Commands}

\begin{document}

\begin{frontmatter}
       
\title{Modeling and simulation of large-scale Systems: a systematic comparison of modeling paradigms}
       
        \author[tug]{G.~Schweiger \corref{cor1}}
        \author[not]{H.~Nilsson}
        \author[kth]{J.~Schoeggl}
        \author[swe]{W.~Birk}
        \author[kfu]{A.~Posch}

        \address[tug]{Graz University of Technology, Graz, Austria}
        \address[not]{University of Nottingham, Nottingham, UK}
        \address[swe]{Lule\aa University of Technology, Lule\aa, Sweden}
        \address[kfu]{University of Graz, Graz, Austria}
        \address[kth]{KTH Royal Institute of Technology, Stockholm, Sweden}

        \cortext[cor1]{Corresponding author. Gerald Schweiger,
        CD Laboratory for Quality Assurance Methodologies for Autonomous Cyber-Physical Systems, Inst. for Software Technology, TU Graz, \textit{E-mail address}: gerald.schweiger@tugraz.at}

\begin{abstract}
A trend across most areas where simulation-driven development is used is the ever increasing size and complexity of the systems under consideration, pushing established methods of modeling and simulation towards their limits.
This paper complements existing surveys on large-scale modeling and simulation of physical systems by conducting expert surveys.
We conducted a two-stage empirical survey in order to investigate research needs, current challenges as well as promising modeling and simulation paradigms.
Furthermore, we applied the analytic hierarchy process method to prioritise the strengths and weakness of different modeling paradigms. 
The results of this study show that experts consider acausal modeling techniques to be suitable for modeling large scale systems, while causal techniques are considered less suitable.  

\end{abstract}


\end{frontmatter}


\input{chapter/Introduction}
\input{chapter/Background}

\input{chapter/Method}
\input{chapter/Results}

\input{chapter/Conclusion}

\section*{Acknowledgement}

The reported research been conducted within the framework of the CD Laboratory for Quality Assurance Methodologies for Autonomous Cyber-Physical Systems at the Institute for Software Engineering at Graz University of Technology. The  financial  support by  the  Austrian  Federal  Ministry  for Digital  and  Economic Affairs and the National Foundation for Research, Technology and Development is gratefully acknowledged.

\section{References}
\bibliography{Mendeley_ModellingParadigms.bib}

\input{chapter/Appendix}

\end{document}

%% file: Commands.tex
\newif\ifcomments
\commentstrue 

%% file: chapter/Introduction.tex
\section{Introduction}

Modeling and simulation based developments have increasingly become established in many research fields and industrial applications. 
Simulation-driven approaches are used to approximate the behaviour of systems and processes in order to improve their efficiency or design new ones.
Furthermore, simplified models on the basis of simulation models are often an integral part of online optimization and control approaches.
A trend across most areas where simulation-driven development is used is the ever increasing size and complexity of the systems under consideration as well as interoperating systems; treating them as islands is not sufficient anymore. 
These trends are pushing established methods of modeling and simulation towards their limits. 

The focus of this paper is modeling and simulation of physical systems that may span multiple physical domains. 
Here, increasing system size and complexity is indeed a pressing concern. Not the least the emergence of so called Cyber-Physical Systems (CPS) poses significant challenges for traditional modeling and simulation techniques. CPS combine computational systems such as microprocessors and communication networks (the ``cyber'' part) with other physical systems \cite{Lee2011a}. Consequently, simulation-driven development of CPS requires combining modeling and simulation techniques for both kinds of system \cite{Kasai2003,Seiger2015}.

There are two broad approaches to modeling and simulation of physical systems: causal and acausal \cite{Fritzson2015}. In causal models, the outputs are explicitly expressed in terms of the inputs; i.e., the direction of information flow is manifest. Concretely, such models are often represented as block diagrams, and ever since the definition of CSSL in 1967, the vast majority of modeling languages for physical modeling has been based on this block-oriented paradigm. However, there are a number of well-known drawbacks to causal modeling, especially at scale \cite{Astrom1998,Cellier1996}. This led to the development of the acausal approach where models essentially are expressed in terms of undirected equations. This makes them much more reusable and composable, addressing some of the challenges of large scale modeling and simulation. 

\subsection{Main contribution}
The main contribution of this paper is to provide a comprehensive discussion about causal and acausal modeling approaches based on a survey of expert opinion. 
We want to stress that the modeling in itself is a process where the intended use of the model often affects the methodology that is used, because certain types of models are not usable with certain tools, like e.g. in control design.
The focus of this paper is on \textit{modeling for simulation}.

To that end, we carried out a study with the participation of 25 world-leading experts in the field of causal and acausal modeling. 
They contributed to the present work by discussing current limitations, future challenges, research needs, as well as strengths and weaknesses of causal and acausal approaches for modeling and simulating large-scale systems.
We applied a two-stage empirical survey and a hybrid method by combining an analysis of Strengths and Weaknesses (SW) with an Analytic Hierarchy Process (AHP), resulting in a SW-AHP analysis.
Using this hybrid method we were able to complement the qualitative results of our study with an aggregated prioritisation of SW-factors. 
To the best of the authors' knowledge, such a research design is novel in the field.


The results and findings of this study can support the efforts of the scientific community to further develop languages, numerical methods and tools for modeling and simulating large-scale systems.
Furthermore, the results can serve as an guide for practitioners in choosing a suitable modeling paradigm for a specific task.

The rest of the paper is organised as follows. Section
\ref{sec:background} surveys relevant literature and gives a more detailed introduction to causal and acausal modeling. Section \ref{sec:method} gives a detailed account of the methods we used to carry out the study. The results of our study are then presented and discussed in section \ref{sec:results}. Finally, we summarise our findings in section \ref{sec:conclusions}.

%% file: chapter/Background.tex
\section{Background}
\label{sec:background}
This section provides an introduction to causal and acausal modeling for physical systems (Section \ref{sec:paradigms}) and a review of the literature concerning the historical perspective as well as critical discussions of both paradigms (Section \ref{sec:background-lr}); this review served as the foundation for the first round of the survey.
\subsection{Paradigms for modeling and simulation of physical systems}
\label{sec:paradigms}


A fundamental distinction can be made between acausal (or non-causal) and causal modeling. Note that ``causal'' here is not used in the sense of \emph{temporal} causality, but rather refers to whether the underlying representation of the modeled system is directed or not, as explained in the following.

In causal modeling, the modeled system is, directly or indirectly, described by a system of ordinary differential equations (ODE) in explicit form; that is, the equations can be viewed as \emph{directed}, making it clear how the unknown quantities are derived from the known ones, hence ``causal´´. In general, such a system of ODEs has the form:
\begin{equation}\label{eq:ode}
\Dot{x}(t) = f(x(t), u(t), p)
\end{equation}
where $x$ is the state (the minimal set of system variables that uniquely determine the future system behavior if the inputs are known),
$t$ is time, $u$ is the control variable, and $p$ are parameters. 
For a directed interpretation, the state derivatives, $\Dot{x}$, are viewed as the unknown quantities, while the remaining ones are considered known.

Causal modeling is also known as explicit modeling, block-oriented modeling, or even imperative modeling as the directed equations can be read as assignment statements. However, the latter does not mean that causal modeling as such entails the use of general imperative state, unlike in imperative \emph{programming}. The term ``block-oriented´´ stems from the fact that such models often are expressed as interconnected blocks, each with designated inputs and outputs. A block can represent a simple mathematical operation (like addition or integration over time), or can be a composite model consisting of interconnected blocks. There is a direct correspondence between such a block diagram and a system of directed ODEs. Signal-flow graphs is another way to represent this kind of models.

In acausal modeling, the modeled system is, directly or indirectly, expressed as a system of differential algebraic equations (DAE) in implicit form. In general, such a system of DAEs has the form:
\begin{equation}\label{eq:dae}
F(\Dot{x}(t), x(t), y(t), u(t), p) = 0
\end{equation}
where $y$ are the algebraic variables and the meaning of $x$, $u$, $t$, $p$ are the same as for the ODE above. As can be seen, there is no longer any manifest directed interpretation. 
Instead, given the known quantities, as determined by a particular usage context, the equations have to be used to solve for the unknowns.

Solving DAEs is, in general, more difficult than solving of explicit ODEs, where the differential variable is a valid system state candidate. 
While the same applies to low-index DAEs, the differentiation index being the number of differentiations needed to convert a DAE into an ODE, the situation is more involved for higher-index problems and other special cases \cite{Pantelides1988}.
However, the DAE formulation means that the modeler can focus more on what to model, rather than on how to model it to, for example, facilitate numerical simulation. Acausal modeling is thus also known as \emph{declarative}, \emph{mathematical}, \emph{physical} or \emph{equation-based}  modeling in the literature \cite{Fritzson2015}. 
Further, the fact that the knowns and unknowns are not given a priori makes the models more reusable; this was a main motivation for the emergence of acausal modeling (see Section \ref{sec:paradigms}).

A further distinction can be made between different solution methods.
State-of-the-art tools based on acausal modeling are based either on methods where equations are globally reduced prior to integration or on the concept of pre-compiled component models. 
Pre-compiled methods facilitate tailored approaches for specific components.
The process of global reduction is described in Cellier and Kofman \cite{Cellier2006} and consists of the steps flattening (transforming the hierarchical structure of a model into a set of DAEs) and causalization (transforming DAEs into ODEs using various algorithms such as matching, sorting and index reduction).

While the terminology varies in the literature, we will, in the following, use the terms causal and acausal in the sense described above: (i) \emph{causal} for the case where the causality or direction is given a priori in a model, and (ii) \emph{acausal} for the case where the causality is not explicitly specified, but inferred when the models are used. Following from this, we will use the term causal modeling language for languages that only support causal models, and acausal or declarative modeling languages for languages that support acausal models.


\subsection{Literature review}
\label{sec:background-lr}




\AA{}str\"{o}m et al.\ \cite{Astrom1998} present a discussion about the historical development of modeling and simulation.
Until the early 90s, engineers have been describing physical systems by linear and nonlinear state space models \cite{ Cellier1996a}. 
Available tools and libraries during that era reflect the emphasis on causal models.
In the 90s, a new modeling paradigm emerged based on several trends and insights: (i) The problems of modeling physical systems based on causal modeling languages \cite{Astrom1998, Cellier1996}; (ii) the demand from user to model and simulate complex multi-domain systems by using object oriented programming languages \cite{ Astrom1998} and  (iii) the demand for reusing models;
in causal modeling languages, a library must contain different types of a mathematical model based on the input-output relation \cite{Cellier1996}.
These drawbacks led to the design and development of equation-based, acausal, object-oriented modeling languages.
In recent years, work has been done in comparing causal and acausal modeling for domain-specific applications. 

Wetter et al.\ \cite{Wetter2016} present a discussion about causal and acausal approaches in building simulation programs.
They conclude that most state-of-the-art building-simulation tools are based on causal paradigms (which are referred to as imperative programming languages). 
While causal paradigms make modeling more difficult and excludes particular powerful methods for simulation and optimization, they identify several advantages for acausal approaches: (i) using equation-based languages allow for a symbolic manipulation of the equations; some numerical simulation methods benefit from the access to these equations; (ii) assessing the equations helps in automatically identifying characteristics of the system that are relevant for simulating stiff systems or hybrid systems; (iii) models based on acausal paradigms are well suited for formulating optimal control problems; using equation-based languages enable tools to convert model equations to a form that is well suited to solve nonlinear optimization problems. 

A similar comparison was conducted by Schweiger et al.  \cite{Schweiger2018c} for energy systems at district scale.
They conclude that acausal modeling is well suited for representing the structure of physical systems and that acausal modeling is convenient for rapid prototyping. 
Furthermore, the conclude that the efficient simulation of large-scale systems is a central consideration for selecting an appropriate modeling language and tool; 
scalability studies are still lacking.
Furthermore, they conclude that standardized models for various applications would be helpful to compare different modeling paradigms, solution strategies and tools.

Wetter and Haugstetter \cite{Wetter2006} reach a similar conclusion. Their study indicates that the model development time in causal languages is five to ten times longer compared with acausal languages. 

Pollok et al.\ \cite{Pollok2019} analyse psychological aspects of acausal modeling approaches.  
They conducted (i) expert interviews, (ii) an experiment based on self-reported timings to analyse the effects of inheritance on ease of understanding, and (iii) an online experiment to analyse the effects of model representations on the performance in modeling tasks. 
They conclude that experienced modelers tend to develop their models from the top-down, while beginners tend to do the opposite. 
The experiments indicate that inheritance significantly increase the time needed to understand a model.
The third experiment indicates that graphical representations reflecting the real-life system structure outperform abstract presentations like block-diagrams for several metrics. 



The literature discusses fundamental and current limitations of modeling languages and tools based on causal paradigms as well as the causes for the development of acausal modeling paradigms. 
A thorough discussion about the advantages and disadvantages of both paradigms as well as their suitability for modeling and simulating large-scale systems are beyond the scope of the existing literature.

%% file: chapter/Method.tex
\section{Method}
\label{sec:method}


To carry out this study, we employed a two stage empirical survey, where we interviewed experts form academia and industry.
Section \ref{sec:method-se} discusses how we selected the group of experts in detail. 

The purpose of the first round was to explore a number of themes emerging from the literature review \ref{sec:background-lr} pertaining to current approaches for modeling and simulation of large-scale physical systems, their limitations, and perceived research needs. 
Most of the questions in the first round were qualitative (i) in order to avoid biased answers, (ii) to introduce the topics in a very broad context and (iii) in order not to miss important perspectives.
The literature review (section \ref{sec:background-lr}) served as the basis for the questionnaire that was put to the experts in round one. 
For the qualitative content analysis of the answers, we followed the framework by Mayring \cite{Mayring2004} for systematic text analysis.



In the second round, most of the questions were quantitative. 
These questions were based on the findings of the first round and the literature review.
For an additional perspective, we also carried out a quantitative investigation of the strengths and weaknesses (SW) of the causal and acausal modeling paradigms applying the Analytic Hierarchy Process (AHP) method in the second round (the method is detailed in Section \ref{sec:SWAHP}).
This method has been applied in a wide range of fields, including renewable energy technologies \cite{Catron2013}, investment behavior \cite{Gottfried2019} or fields similar to those of the current study (e.g. co-simulation techniques \cite{Schweiger2019a}).

\subsection{SW-AHP}
\label{sec:methSWAHP}
In the second step, we conducted an Analytic Hierarchy Process (AHP) analysis, to identify the relative importance of Strengths and Weakness (SW) factors for causal and acausal modeling. 
With Strengths and Weaknesses (SW) we refer to positive and negative factors that we were able to identify in the course of the first round of interviews and the literature review.
We limited the number of SW factors to four per group (i.e. four per Strengths/Weaknesses). This was done to ensure sound analysis and to enable experts to complete the questionnaires within an acceptable period of time. 

In order to verify that the most important factors for each group were identified for SW-AHP analysis, we added open questions, where experts could list further factors 
In the course of the AHP, experts undertook a pair-wise assessment of the factors for the respective group;
 i.e. to state which Strengths/Weakness factor for each pair is more important, and how much more important. 
For all the comparisons, we applied the nine step scale suggested by Saaty \cite{Saaty1986}. It ranges from $9:1$ (meaning that some factor $a$ under consideration is much more important than a factor $b$ under consideration) to $1:9$ (meaning that factor $b$ is much more important than factor $a$). 
The even numbers are omitted as intermediate steps. 
The center of the scale is $1:$1; this means that the respective factors were regarded equally important. 

Because we had four factors in each SW group,
the respondents were asked to make 12 pairwise comparisons in total, six for each of the SW fields. Each pairwise comparison followed the logic shown below for the comparison of factor $a$ and $b$. For example, where respondent $i$ compares the factors $a$ and $b$, and judges factor $a$ to be more important than factor $b$, then the weighting of factor $a$ with respect to $b$ by respondent $i$, written as $w(a)_i/w(b)_i$, is an odd number between $3$ and $9$, depending on $i$'s judgment of $a$'s relative importance with respect to $b$. Where factor $a$ is judged to be of the same importance as factor $b$, the value $w(a)_i/w(b)_i$ is equal to $1$. Where factor $a$ is judged to be less important than factor $b$, the value $w(a)_i/w(b)_i$ is the reciprocal value of the odd numbers between $3$ and $9$.


We then calculated the average figure for the results of all pairwise comparisons. Here, we normalized the average scores of each comparison between two factors so that the less important factor always received a score of $1$, and the more important factor a score within the possible range from $1$ to $9$. 

\subsection{Selection of experts}
\label{sec:method-se}

The selection of experts who participated in the survey was based on a so called Knowledge Resource Nomination Worksheet (KRNW) \cite{Okoli2004}, which was developed by Delbecq et al. \cite{Delbecq1975GroupProcesses} to select experts within a nominal group technique. 

It involves the following five stages: (1) preparing the KRNW by selecting experts from industry and academia; (2) populating the KRNW; we populated both group of experts based on a keyword-based literature study (i.e., if they had (co-)authored a publication).
Additionally, the category of industrial experts was further populated based on the experience of the authors and based on experience with consulting practitioners. (3) Nominating further experts based on the feedback from the initial experts. 
(4) Ranking of the experts by number of publications and citations \footnote{www.scopus.com}.
Scopus is the world's largest scientific database for peer-reviewed literature \cite{Gonzalez2010}.
(5) Inviting the experts to participate via an online questionnaire. 

32 experts were contacted to complete the questionnaire. 25 experts completed the questionnaires; the response rate was 78\%. A total of 8 experts from industry and 14 from the university sector took part; 3 did not provide any information. 
It should be stressed that the number of experts participating in this survey is consistent with the guidance given in the literature. \cite{Ziglio1996,Clayton1997,Ludwig1997}.

The experts were asked to estimate their level of expertise in causal and acausal modeling on a four-level Likert scale (4 = high, ..., 1 = none). Results show that more than 80\% of the surveyed experts have a high or moderate level of expertise in both modeling approaches (causal modeling: Interpolated Median (IM; see Section \label{sec:method-pr}) = 3.2; acausal modeling: IM = 3.5).
Thus, biased answers could be avoided. 
Industrial experts who participated in the second round can be assigned to the following fields: Software development: 3, Automotive: 1, large-scale simulation development: 1, test and measurement: 1, computer science: 1, energy related applications: 1. 

Academic experts who participated in the second round can be assigned to the following research fields: Software development: 4, energy related applications: 3, Computer Science: 2, Physiology: 2; Automotive: 1, Maritime: 1,  Algorithm development: 1. 

The experts thus had fairly heterogeneous backgrounds, and we can conclude that the size of our expert group is aligned with the accepted recommendations.

\subsection{Presentation of the results}
\label{sec:method-pr}
Hallowell and Gambatese \cite{Hallowell2010QualitativeResearch} suggest that the median value is better suited for presenting results than the mean value;
the median value is less affected  by extreme values of outliers.
Furthermore, Sachs \cite{Sachs1997} suggests that the interpolated median is more accurate than the median value.
It gives a measure within the upper bound and lower bound of the median, in the direction that the data is more heavily weighted.

In this paper the most important results are presented in a bar chart in the Appendix (see Section \ref{sec:appendix}).
Furthermore, we analyse the responses in terms of their median
(M), interpolated median (IM ), and the average (A).
This should ensure a transparent presentation of the findings.

The $\mathit{IM}$ is calculated as:
\begin{equation}\label{eq:first}
\mathit{IM} =
\begin{cases} 
      M & \text{if $n_2=0$}, \\
      M-0.5+ \frac{0.5 \cdot N-n_1}{n_2} & \text{if $n_2\neq 0$}
   \end{cases}
\end{equation}

$N$ is the number of answers to a specific question;
$n_1$ is the number of answers strictly less than $M$; $n_2$ is the number of answers equal to $M$.


\subsection{Threats to validity} \label{threats}

There is no generally accepted measure that allows a transparent comparison of the productivity and the impact of a researchers' work.
This is especially the case for the comparison between disciplines and research fields.

The academic experts for this survey were selected exclusively on the basis of the number of publication.
The authors are aware that this is a threat to validity;
nevertheless, the authors assume that this is the most transparent selection process.
Industrial experts were also selected based on the number of publications as well as the authors' experience with consulting practitioners. 
This ensured that industrial experts who do not publish their work in scientific journals nevertheless took part in the survey.
This selection process could be considered as unrepresentative.
The answers of the experts show, however, that they have great expertise in the field of modeling and simulation of large-scale systems.

%% file: chapter/Results.tex
\section{Results and Discussion}
\label{sec:results}
This chapter presents the main results of the empirical survey and the SW-AHP analysis. It is structured as follows. In Section \ref{sec:utilization}, the utilization of causal and acausal languages and tools is presented. 
Section \ref{sec:large_scale_systems} presents the expert assessment of whether it would be helpful to define a notion of large-scale models to facilitate comparison of different modeling paradigms, solution strategies and tools, and how feasible it would be in practice to define such a notion. 
Section \ref{sec:comp_causal_acausal} then presents the expert assessment regarding causal and acausal modeling and simulation of large scale systems. 
This includes the experts' assessment of main limitations, most important levers for simulating large-scale systems as well as the potential of specific approaches and techniques for improving the efficiency of simulation of large-scale models, both in causal and an acausal settings. 
Section \ref{sec:SWAHP}, finally, presents the results of the SW-AHP analysis.

\subsection{Utilization of causal and acausal languages } 
\label{sec:utilization}
The utilization of tools based on causal and acausal modeling approaches is investigated using a keyword analysis on Scopus.
The candidates were identified based on the knowledge of the authors and the relevant literature such as \cite{Astrom1998, Fritzson2015,Schweiger2018c, Pollok2019}.
We have only considered those languages that have listed more than 50 publications with the respective keyword on Scopus since 2000.
It should be emphasized that some acausal languages also support causal modeling. 
However, the authors’ experience is that these are not yet competitive in the field of causal modeling with strictly causal approaches. 
Therefore, both paradigms are considered separately in this work.

\begin{figure}[h!]
\centering
\includegraphics[width=1\textwidth]{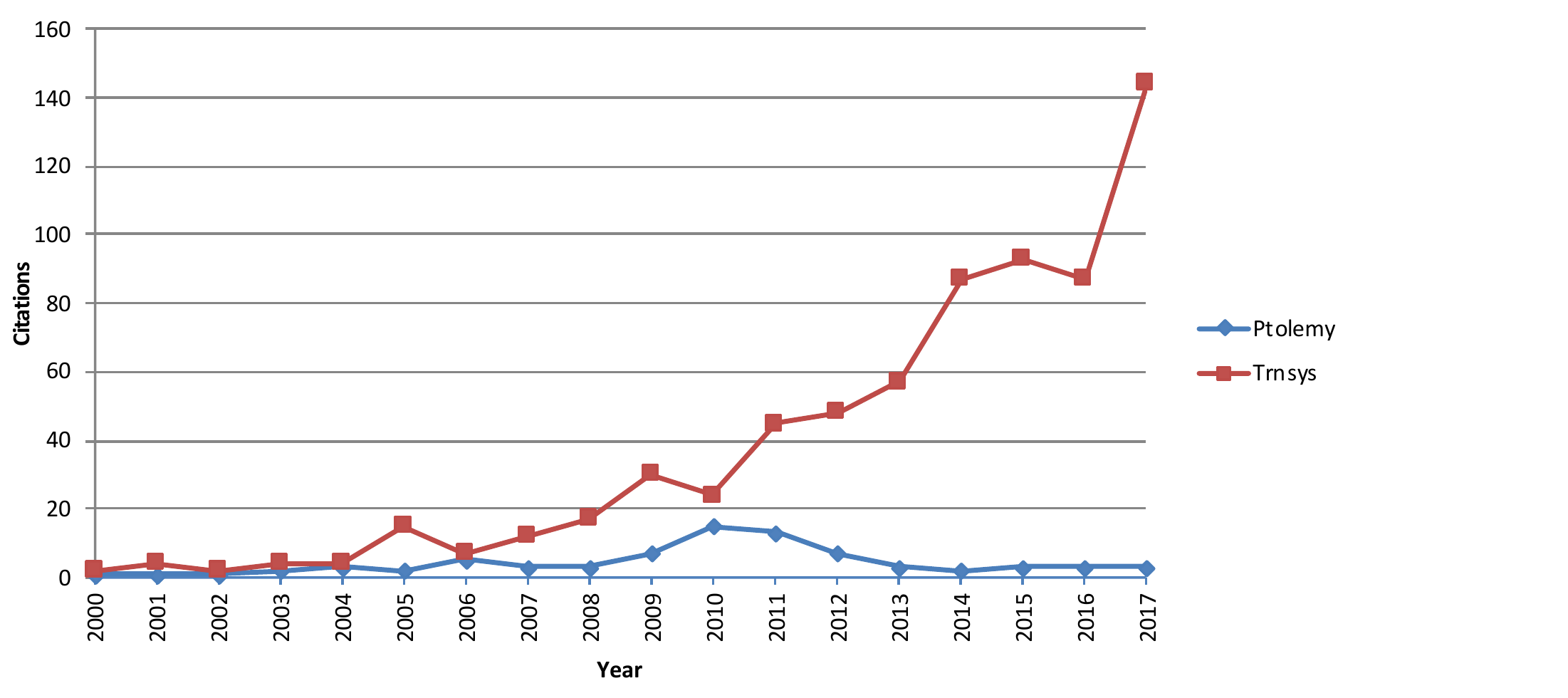}
\caption{Utilization of causal modeling languages and tools}
\label{fig:causal_ut}
\end{figure}

\begin{figure}[h!]
\centering
\includegraphics[width=1\textwidth]{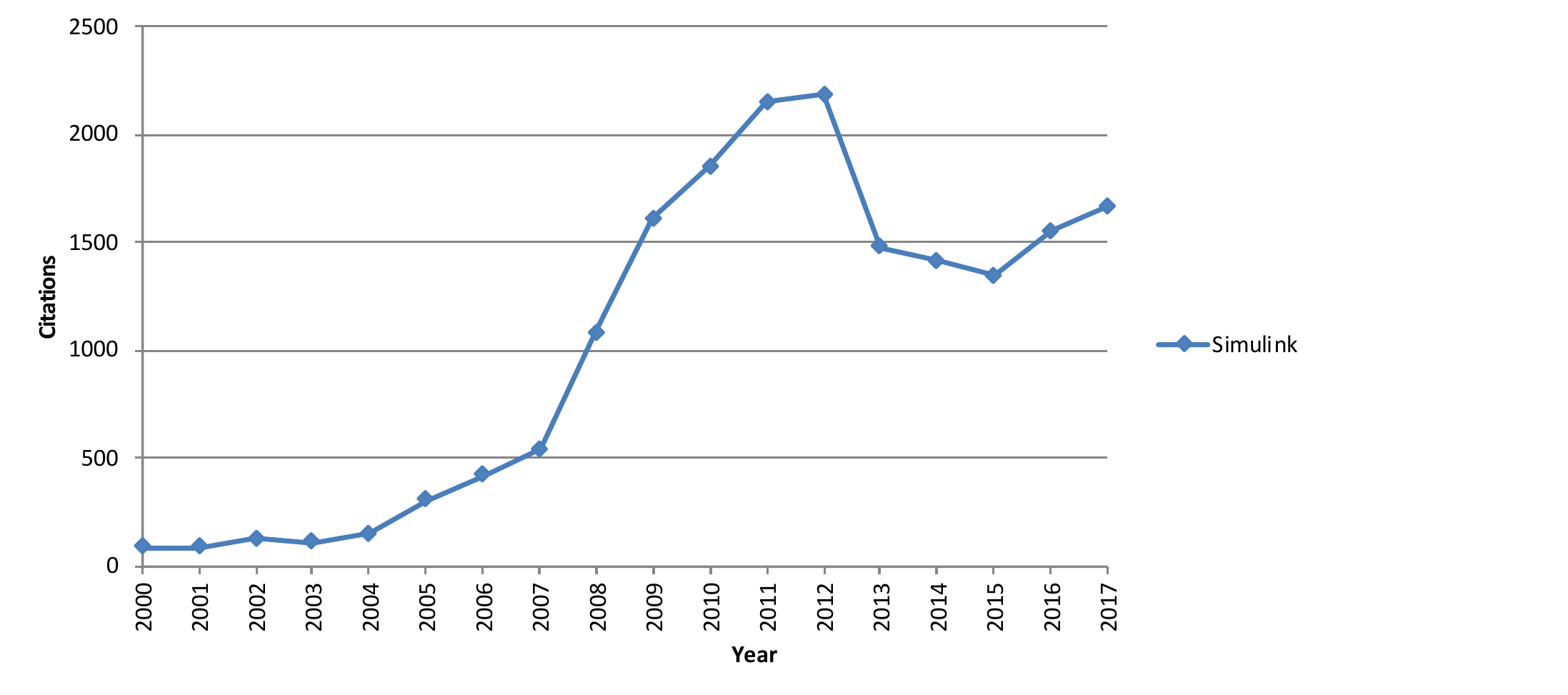}
\caption{Utilization of Simulink}
\label{fig:Simulink}
\end{figure}

Overall, Simulink is by far the most widespread (see Figure \cref{fig:Simulink}).
The number of publications per year was over 2000 in 2011 and fell below 1500 in 2015.
As can be seen from \cref{fig:causal_ut}, there is an increasing utilization of TRNSYS since 2006 (approx. 140 publications in 2017).
However, publications with the keyword Simulink are more than 10 times more frequent in 2017 than those with the keyword TRNSYS.

\begin{figure}[h!]
\centering
\includegraphics[width=1\textwidth]{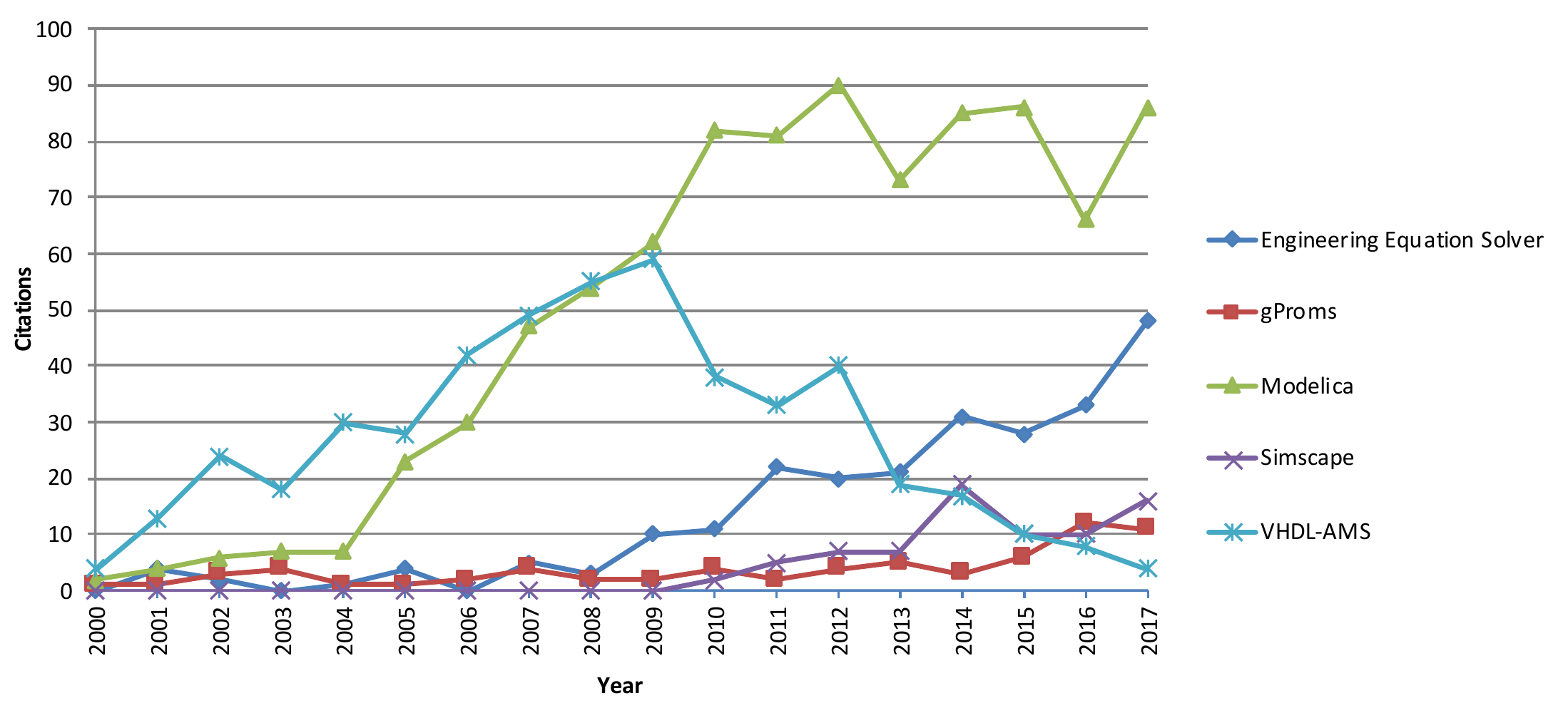}
\caption{Scopus keyword search; total number of publications for various tools based on acausal modeling approaches (left axis); number of Modelica publications for different domains (right axis). The amount of publications based on Modelica is roughly constant over the recent years}

\label{fig:acausal_ut}
\end{figure}

As can be seen from \cref{fig:acausal_ut}, the most widely used language for acausal modeling is Modelica. 
It should be noted that Modelica is a modeling language supported by a number of different modeling and simulation tools.
While there was an increase in the number of publications between 2004 and 2010, the number has stagnated since then and is currently between 70 and 80 per year.
These results do not confirm that acausal approaches have a growing user base; this was indicated e.g. in \cite{Pollok2019} for Modelica.
Publications with the keyword Simulink are about 20 times more frequent in 2017 than those with the keyword Modelica.
No other causal or acausal language has more than 50 publications with the respective keyword in 2017.

\subsection{Standardizing a notion of large scale systems}
\label{sec:large_scale_systems}

The experts were asked if a standardized notion of what it means for a model to be large-scale would be helpful, and how feasible it would be to define such a notion. The respondents could choose in a five-level Likert scale (5: Strongly Agree; 4: Agree; 3: Neither Agree nor Disagree; 2: Disagree; 1: Strongly Disagree). In the first round, experts were asked, in an open-ended question, about a suitable measure to define large-scale models in terms of metrics and properties of the model. We analyse the responses in terms of their median ($M$), interpolated median ($\mathit{IM}$), and, sometimes, the average ($A$); see section \ref{sec:method-pr}. The responses show that the expert view is that having a standardized notion of what it means for a model to be large-scale would be both helpful ($M= 5$, $\mathit{IM} = 4.7$) and feasible ($M = 4$, $\mathit{IM} = 4.1$). 

However, the answers regarding possible metrics and properties of the models varied widely, ranging from (i) number of states, nonlinear degree, size of loops, dynamism to (ii) number of connected sub-systems, to (iii) a metric should probably be constructed from the product of the computational complexity of determining the rate vector and the computational complexity of integrating it. Further, many of the experts were skeptical regarding the possibility of defining such a standardized notion at a sufficiently high level of abstraction so as to be applicable across a wide range of kinds of systems.

As a representative example of ideas put forward in the first round, one expert held that it would be useful and possible to define several classes of large scale problems for \emph{specific} domains allowing the comparison of causal and acausal modeling approaches and solution algorithms. These classes would be defined on an abstract level (e.g., continuous/hybrid, highly nonlinear/linear, sparse/dense) and then translated in domain specific problems. As examples, the expert suggested two applications in the area of Smart Grids: (i) A continuous problem of $n$ buildings implemented as $RC$ models with electric cooling, connected to a power distribution grid and $n$ producers; (ii) Converting the previous problem into hybrid problem by including state and time events and linearize all sub-models (buildings, HVAC systems, etc.). The equations for all sub-models as well as modeling assumptions and boundary conditions should be presented in a transparent manner.

Some other concerns, besides how widely applicable a single standardized notion of large scale system could be, were also raised. One expert pointed to a possible bias in defining suitable characterizations: if these were to be developed by tool vendors or certain modeling language experts, the characterizations may be developed in a way that showcases their preferred approach. This expert therefore argued that it would be better if the characterizations were done by domain experts, who can be expected to be more neutral, even if domain experts in many cases also would be experts in certain languages and tools. Another point this expert raised was that the implementation should be done in open competition as the actual implementation in certain languages and tools can have a significant impact on the performance of the simulation. 

\subsection{Comparison of causal and acausal modeling and simulation}
\label{sec:comp_causal_acausal}
Experts were asked to assess the suitability of causal and acausal approaches to model and simulate large-scale systems;
table \ref{tab:LS} summarizes the responses.

\begin{table}[htbp]
\tiny
  \centering
  \caption{Experts assessments: Suitability of causal and acausal approaches to model and simulate large-scale systems. Scale: 7: Entirely Agree; 6: Mostly Agree; 5: Somewhat Agree; 4: Neither agree nor disagree; 3: Somewhat disagree; 2: Mostly disagree; 1: Entirely disagree. Headings: $A$: Average, $M$: Median, $IM$: Interpolated Median.}
    \begin{tabularx}{1\linewidth}{
        X
        >{\color{black}}R{0.8cm}
        >{\color{black}}R{0.8cm}
        >{\color{red}}R{0.8cm}
        <{\color{black}}
    }
    Question & \thead{$A$} & \thead{$M$} & \thead{$IM$} \\
    \midrule
    Causal modeling techniques are suitable for modeling large scale systems  & 3.5   & 3.0   & 2.0 \\
    Causal models are suitable for simulating large scale models  & 4.8   & 6.0   & 5.8 \\
    Acausal models are suitable for simulating large scale models & 5.9   & 6.0   & 6.1 \\
    Acausal modeling techniques are suitable for modeling large scale systems & 6.5   & 7.0   & 6.8 \\
    \bottomrule
    \end{tabularx}%
  \label{tab:LS}%
\end{table}%

In terms of suitability for modeling large systems, experts see a clear difference between causal and acausal modeling. There was a strong agreement ($M = 7$, $\mathit{IM} = 6.8$) that acausal modeling techniques are suitable for modeling large scale systems, 
while causal techniques are considered less suitable ($M = 3$, $\mathit{IM} = 2.0$): more than 50 \% of the experts disagreed or disagreed entirely with the statement that causal modeling techniques are suitable for modeling large
scale systems. 
When it comes to simulation, both approaches were seen to more or less equally suitable, reflecting the fact that it is possible to construct large models using either approach, and once the hard work of constructing a model suitable for simulation has been done, it does not matter so much how one got there.
One expert emphasised that causal modeling is a proper subset of acausal modeling which can support both causal and acausal formalisms; 
Causal modeling cannot replace acausal modeling in terms of functionality; both are required, but both can be provided via an acausal approach.

An interesting case are emerging CPS that combine cyber systems with physical systems.
One expert mentioned that it may be easier to use acausal modeling for describing physical processes, especially when the causality is not well defined statically. However, when designing systems that have a cyber component, the designer usually has causality relationships on the basis of information flow in mind; as an example, the expert mentioned ``the controller actuates via actuators, influences some physical process, which is then sensed by the sensors and given as inputs to the controller. At a high level view, I can draw arrows between my components (controller, actuator, plant, sensor), even though, at a lower level, I would use acausal modeling for the plant (the physics), and causal modeling for the controller.'' This is a typical closed loop control context, where feedback loops are used to achieve predefined performance characteristics of the process in question. The overall behavior of the system is therefore reflected by the combined model comprising both the cyber and physical part.


\subsubsection{Acausal modeling and simulation} 
\label{sec:acausal_modelling}
\underline{Potential}:
Experts were asked about the potential of acausal approaches for large-scale modeling and simulation.
Several experts mentioned that the modeling process itself within acausal approaches is the greatest strength. 
Several experts pointed out that hybrid acausal modeling touches an enormous range of potential problems in engineering.
While one expert stressed that state of the art acausal tools can provide performance that is comparable to special purpose tools, another expert pointed out that in case symbolic processing cause scalability problems, parts of the acausal model could be converted to causal models. 
Hence, acausal models should have all the benefits of causal models, plus additional benefits from acausal appraoches.
Furthermore, experts stressed that the integration with other modeling paradigms such as agent-based or  FEM-like/CFD-like simulations into systems simulations shows great potential.

Due to many detailed answers in the first round concerning specific concepts and techniques to enable or improve the efficient simulation of large-scale models, experts were asked in the second round to rate the potential of these concepts and techniques based on acausal approaches. 
The results are summarized in table \ref{tab:AcausalImpr}. We discuss the concepts and techniques the experts consider having high potential ($\mathit{IM} \geq 3.5$ ) in more detail below.

\begin{table}[htbp]
\tiny
  \centering
  \caption{The potential of specific concepts and techniques for enabling or improving efficient simulation of large scale systems. Scale: 4: High Potential; 3: Moderate Potential; 2: Low Potential; 1: No Potential.}
    \begin{tabularx}{1\linewidth}{
        X
        >{\color{black}}R{0.8cm}
        >{\color{black}}R{0.8cm}
        >{\color{red}}R{0.8cm}
        <{\color{black}}
    }
    Concept/Technique & \thead{$A$} & \thead{$M$} & \thead{$\mathit{IM}$} \\
    \midrule
    Linear system solvers & 2.1   & 2.0   & 2.1 \\
    Exploiting repetitive structures & 3.1   & 3.0   & 3.1 \\
    Multi-Rate Event Handling & 3.3   & 3.0   & 3.2 \\
    Index reduction methods & 3.0   & 3.0   & 3.3 \\
    Multi-Rate Algorithms & 3.4   & 3.0   & 3.3 \\
    DAE Solvers & 3.4   & 4.0   & 3.6 \\
    Parallel computing & 3.6   & 4.0   & 3.7 \\
    QSS algorithms & 3.6   & 4.0   & 3.7 \\

    \bottomrule
    \end{tabularx}%
  \label{tab:AcausalImpr}%
\end{table}%

According to the experts, the greatest potential for improving simulation of large scale systems is offered by quantized state systems (QSS) methods and parallel computing ($\mathit{IM} = 3.7$ in both cases), followed by DAE solvers ($\mathit{IM} = 3.6$). In all three cases, around 60 \% of the experts considered the potential as high.

QSS are an alternative class of algorithms for numerical integration based on discretizing the state values, while keeping time continuous \cite{Cellier2006}. This algorithm turn ODE systems into Discrete Event Systems (DEVS). State of the art QSS algorithms include second- and third-order methods as well as a linearly implicit method that is suitable for stiff systems \cite{Bergero2012}. 


Comparisons of QSS with state of the art solver based on discretion over time show promising results. For example, Floros et al. \cite{Floros2011} compared a QSS third-order algorithm and a DASSL solver in OpenModelica. The QSS method could reduce the CPU time by a factor of 20. Bergeo et al. obtained similar results \cite{Bergero2012}: the stiff LIQSS solver could reduce the CPU times by a factor of 40 compared to DASSL.

In recent work, new algorithms were proposed to transform DAEs from a high level modeling language to a special index one form  that can be solved with state-of-the-art Index-1 DAE integrators. 
This algorithms keep the sparsity of the model equations; no equations systems are solved in the transformation system \cite{Otter2017, Braun2017}.

Recent studies \cite{Casella2015} have emphasized that in many large-scale models such as smart grids, subsystems only interact with a few nearby subsystems. 
In these systems, the derivative only depends on the value of a few neighbouring states. 
Thus, the Jacobian has only a few non-zero elements meaning it is sparse. 

The inefficiency of using dense solvers for sparse systems has been discussed by Casella \cite{Casella2015} and Braun et al. \cite{Braun2017} among others. 
Sparse ODE and DAE solvers are available. 
While ODE solvers can be both explicit and implicit, DAE solvers are always an implicit, having to solve a non-linear system in each iteration. 
The Jacobian gets bigger in DAE methods, since the integration methods needs to solve for $\dot(x)$, $x$ and $y$ instead only for $x$. 
Braun et al. \cite{Braun2017} show that sparse solvers outperform traditional dense solver for many applications. 
However, no general conclusion could be drawn for the performance comparison of sparse DAE solvers and sparse ODE solvers: the results depend on the specific problem.

A key question when using parallel computing is how to take advantage of the latest advances in multicore technology that lead to enormous computing power, in a user-friendly way (partitioning needs to be made automatically).
In recent years, CPUs as well as GPUs are used in a variety of applications and both have different characteristics and strengths \cite{Mittal2015}.
Elmqvist et al. \cite{Elmqvist2014} developed a method to automatically parallelize model equations in Modelica;
results show that a speed-up of 3.4 times has been achieved using 4 cores/8 threads.
Elmqvist et al. \cite{ Elmqvist2015a} presented a Modelica code generation for GPU cores.
This is in particular important when the model has regular structure, for example, discretization of PDEs, where each cell can then be characterized by a function call.
This evaluation can be made in parallel on GPU cores.
The results of a prototype implementation showed that the use of GPU speeds up the simulations by a factor of five.
Carl and Biswas \cite {Carl2016} investigated methods to partition and simulate differential equation-based models of cyber-physical systems (using a subset of Modelica) using multiple threads on multi-core CPUs;
The main result of this paper are recommendations on (i) the size of the model to be parallelized, (ii) the number of threads to use for the simulation,(iii) the memory management scheme to use, (iv) how to divide the computational work between the threads, and (v) software optimizations to implement
Casella \cite{ Casella2015} pointed out that parallel computation strategies should be combined with appropriate techniques exploiting sparsity and locality of large-scale models, to achieve significant performance improvement. 
Furthermore, Casella mentioned that performance tests in recent literature are carried out with low-end computers with few parallel cores. 


\underline{Main limitations}:
Experts were asked about the main limitations of state of the art languages and simulation environments based on acausal modeling.
Several experts pointed out that a major limitation is the limited debugging functionality: it is difficult to find model errors; a main reason is that only global error feedback is provided.

Further limitations were mentioned: (i) huge experience is required to speed the simulation up, (ii) languages are inflexible and lacking clear semantics, (iii) some tools (e.g. Modelica tools) use one time step for all equations which is not good for scalability of stiff systems and (iv) there are missing style guides to model large-scale systems.
Several experts have commented that these points are not a deficiency of using an acausal language, but rather a deficiency of the current tools.

Based on comments of experts in the first round, experts were asked in the second round about differences in solution strategies: (a) pre-compiled component models vs. (b) globally reducing equations prior to integration.

The answers from the experts clearly showed that the pre-compiled component models has many advantages and few disadvantages compared with globally reducing equations prior to integration.
The following advantages of pre-compiled component models were mentioned: (i) although it is harder to describe what a compilation unit is in modeling (template of equations) vs. software (functions),
separate compilation allows the cost incurred by the computational
complexity of compilation algorithms to be managed;
(ii) it allows for IP protection;
(iii) less start-up time for simulation;
(iv) it helps for non-linear sub-components such as agents;
(v) faster simulations.

The following disadvantages were mentioned: (i) it is not insightful; (ii) implementing separate compilation is more difficult and requires a way of representing the compiled models.

Experts mentioned the following advantages for methods where equations are globally reduced prior to integration: (i) it is more flexible;
(ii) increased re-use of components; 
(iii) it may be faster on long running simulations
(iv) the later equations are turned into simulation code, the more flexible the simulation language.
(v) it is more insightful.

\underline{Lever}: 
Experts were asked about the most important lever for simulating large-scale systems based on acausal appraoches more efficiently.
Experts gave the following answers regarding the most important levers in the near future (1-5 years): (i) Improved debugging resulting in easier readable error messages;
(ii) open development process that empower people to contribute;
(iii) improved symbolic and numerical Methods.

Experts gave the following answers to the most important levers in the next 5-10 years: 
(i) better hardware;
(ii) parallel computing;
(iii) One expert mentioned that compared to other technology eco-systems, state of the art approaches are missing a flattened representation of hybrid acausal systems.  Acausal languages are analogous to languages like C++/Java etc. 
They include lots of high-level language features that make it easier for (human) developers to describe and compose problems. 
The expert concluded that ``what is missing is something analogous to LLVM/WebAssembly/Object code that represents the simplest possible representation of the problem (devoid of inheritance, hierarchies, etc). 
I think this is an important strategic development because it could help to 
\textit{divide and conquer} methods development from tool development. 
This would allow more of a marketplace for \textit{solvers} that is independent of the marketplace for \textit{authoring environments}.
Right now, those are too tightly coupled in the modeling space.''

\subsubsection{Causal modeling and simulation} 
\label{sec:causal_modelling}
\underline{Potential}:
Experts were asked about the potential of causal approaches for large-scale modeling and simulation.
Several experts mentioned that causal approaches are well suited for specific well known problems, mainly in the control domain.
Furthermore, experts mentioned that a causal approaches are wide spread. 

\underline{Main limitations}:
Experts were asked about the main limitations of state of the art languages and simulation environments based on causal modeling.
Several experts pointed out that a major limitation for causal approaches is the need for causalization of the model's equations; it is seen as a main limitation to describe something in a causal manner what is generally an inherently acausal.
Furthermore, experts assess causal approaches as inflexible, difficult to read, write and extend models.


\underline{Lever}:
Experts were asked about the most important lever for simulating large-scale systems based on causal approaches more efficiently.
Generally, many experts mentioned that causal methods are already very well developed and there is hardly any room for improvement.
Experts gave the following answers regarding the most important levers in the near future (1-5 years): (i) multi-rate solvers (ii) explicitly integrate models and solvers together for specific cases.


\subsection{Analytic Hierarchy Process of Strengths and Weaknesses}
\label{sec:SWAHP}
As outlined in \cref{sec:methSWAHP}, in the last research step, an Analytic Hiearchy Process (AHP) was conducted to assess the relative importance of the SW-factors (Strengths and Weaknesses) of the two modeling paradigms that were compiled in the preceding step. 
The SW-AHP analysis is shown in \Cref{table:acausal} for acausal modeling and in \Cref{table:causal} for causal modeling. 
The Strenths factor \textit{S3} of acausal modeling ``Access model knowledge that is stored in equations'' received the highest relative priority of 0.37, followed by \textit{S1} ``Well suited for rapid prototyping'' (priority = 0.27). The third and fourth most important Strengths factors are \textit{S2} ``Easy to read and interpret models'' (priority = 0.23) and \textit{S4} ``Useful in education'' (priority, 0.13). This result indicates that the functionality/performance-related Strengths of acausal modeling paradigms dominate over the accessibility of the models and their use for educational purposes.

A similar distribution of priorities among the four Weakness factors can be observed. 
Factor \textit{W2} ``Debugging of acausal models'' received the highest priority (0.40) in this group, rendering it the central Weakness of acausal modeling paradigms. The ``Small community in industry and academia'' \textit{W4} was found to be the second-most important Weakness (priority = 0.23), followed by the factor \textit{W1} ``Limited suitability for efficiently simulating large scale systems'' (priority = 0.21). The Weakness factor \textit{W3} ``Lack of education and material'' was considered the least important. 


\begin{table}[]
\resizebox{\textwidth}{!}{%
\begin{tabular}{lll}
\hline
\textbf{Strengths}  &                                                                    & Priority \\ \hline
\textit{S3}         & Access model knowledge that is stored in equations                 & 0.37           \\
\textit{S1}         & Well suited for rapid prototyping                                  & 0.27           \\
\textit{S2}         & Easy to read and interpret models                                  & 0.23           \\
\textit{S4}         & Useful in education                                                & 0.13           \\ \hline
\textbf{Weaknesses} & \textbf{}                                                          & Priority \\ \hline
\textit{W2}         & Debugging of acausal models                                        & 0.40           \\
\textit{W4}         & Small community in industry and academia                           & 0.23           \\
\textit{W1}         & Limited suitability for efficiently simulating large scale systems & 0.21           \\
\textit{W3}         & Lack of education and material                                     & 0.15           
       
\end{tabular}%
}
\caption{Result of the SW-AHP for acausal modeling. SW = Strenghts and Weaknesses. AHP = Analytic Hierarchy Process. Priority refers to the relative priority of a single factor within its respective group; i.e. Strengths or Weaknesses. In each group, the priorities add up to 1.}
\label{table:acausal}  
\end{table}

As it is illustrated in \cref{table:causal}, The SW-AHP for causal modeling paradigms revealed, that the Strengths with the highest priority is factor \textit{S4} ``Well suited for the design and implementation of control schemes'' with a priority of 0.42. Factor \textit{S4} also received the highest priority of all assess SW-factors and it was considered almost twice as important as factor \textit{S3} ``Widely used in industry and academia'' (priority 0.22) and factor \textit{S1} ``Well suited for efficiently simulating large scale systems (models are close to the solution algorithm)'' (priority = 0.21). The ``Commercial tool support'', factor \textit{S2} ranked least with a priority of 0.16. 

In the group of Weaknesses of causal modeling, factor \textit{W1} ``Limited suitability for rapid prototyping (difficult to extend, adopt and reuse models)'' ranked first (priority 0.40). With a priority of 0.23, the Weakness factor \textit{W2} ``Difficult to read and interpret models'' ranked second, followed by the ``Causal modeling requires more modeling knowledge than acausal modeling'' \textit{W4}. Factor \textit{W3} ``Causal modeling takes more time than acausal modeling'', was considered as the least important weakness. This result illustrates that that the expert panel considered the higher demand in knowledge and time that causal modeling has in relation to acausal modeling, not as critical as the paradigm's limitations regarding rapid prototyping and model interpretation. When comparing the SW-AHP results for the two paradigms the pivotal role of a model's ``suitability for rapid prototyping'' becomes apparent. It is he second-most important Strength of acausal, factor \textit{S1}, and the limitations regarding it, most important Weakness of causal modeling paradigms, factor \textit{W1}. 

\begin{table}[]
\resizebox{\textwidth}{!}{%
\begin{tabular}{lll}
\hline
\textbf{Strengths}  &                                                                    & Priority \\ \hline
\textit{S4}         & Well suited for the design and implementation of control schemes  & 0.42           \\ 
\textit{S3}         & Widely used in industry and academia                 & 0.22           \\
\textit{S1}         & Well suited for efficiently simulating large scale systems (models are close to the solution algorithm  & 0.21           \\
\textit{S2}         & Commercial tool support       & 0.16           \\ \hline
\textbf{Weaknesses} & \textbf{}                                                          & Priority \\ \hline
\textit{W1}         & Limited suitability for rapid prototyping
(difficult to extend, adopt and reuse models) & 0.40           \\
\textit{W2}         & Difficult to read and interpret models                                        & 0.26           \\
\textit{W4}         & Causal modeling requires more modeling knowledge than acausal modeling                           & 0.20 \\
\textit{W3}         & Causal modeling takes more time than acausal modeling                                     & 0.14

\end{tabular}%
}
\caption{Result SW-AHP for causal modeling. SW = Strenghts and Weaknesses. AHP = Analytic Hierarchy Process. Priority refers to the relative priority of a single factor within its respective group; i.e. Strengths or Weaknesses. In each group, the priorities add up to 1.}
\label{table:causal} 
\end{table}

%% file: chapter/Conclusion.tex
\section{Conclusion}
\label{sec:conclusions}
A trend across simulation-driven development is the ever increasing size and complexity of the systems under consideration, as well as the increasing prevalence of interoperating systems. This pushes established methods of modeling and simulation towards their limits.
This paper presents a two stage empirical survey in which experts assessed current challenges and remaining research needs of modeling paradigms for modeling and simulation of large scale physical systems that may span multiple physical domains.

The main findings from this survey are: 
\begin{compactitem}
\item Experts consider a standardized notion of what it means for a model to be large-scale to be helpful and feasible. 
The definition of this problem should be transparent and unbiased.
\item Experts consider acausal modeling techniques to be suitable for modeling large scale systems,  while causal techniques are considered less suitable; when it comes to simulation, both approaches are seen to be similarly suitable.
\item Experts see the greatest potential for improving the simulation of large-scale systems based on acausal modeling techniques in QSS methods and parallel computing, followed by DAE solvers. Furthermore, experts consider an open development process as an important lever for acausal approaches, while improvements to debugging methods for acuasal modeling is important to enable such techniques to reach their full potential.
\item Many experts mentioned that causal methods are already very well developed and there is hardly any room for improvement.
\item The strengths factor ``Access model knowledge that is stored in equations'' received the highest relative priority in the SW-AHP analysis for acausal modeling techniques. The factor ``Debugging of acausal models'' received the highest priority in the weaknesses group.
\item The strengths factor ``Well suited for the design and implementation of control schemes'' received the highest relative priority in the SW-AHP analysis for causal modeling techniques. The factor ``Limited suitability for rapid prototyping'' received the highest priority in the weaknesses group.
\end{compactitem} \leavevmode

Causal modeling is in principle a special case of acausal modeling and the latter could therefore be used as the general approach. Nevertheless, causal modeling is in some contexts more efficient and therefore the preferred approach.
Harmonization of the modeling formalisms would therefore increase the modeling efficiency and should be investigated.

When comparing the SW-AHP results for the two paradigms the pivotal role of ``Suitability
for rapid prototyping'' becomes apparent. It is the second-most important strength of acausal modeling, and the limitations regarding it the most important weakness of causal modeling.

We hope that the results of this study will stimulate further studies in this field (e.g. in-depth analysis of modeling requirements for CPS; modeling and optimization of large-scale systems).

%% file: chapter/Appendix.tex
\section{Appendix}
\label{sec:appendix}

\begin{figure}[h!]
\centering
\includegraphics[width=1\textwidth]{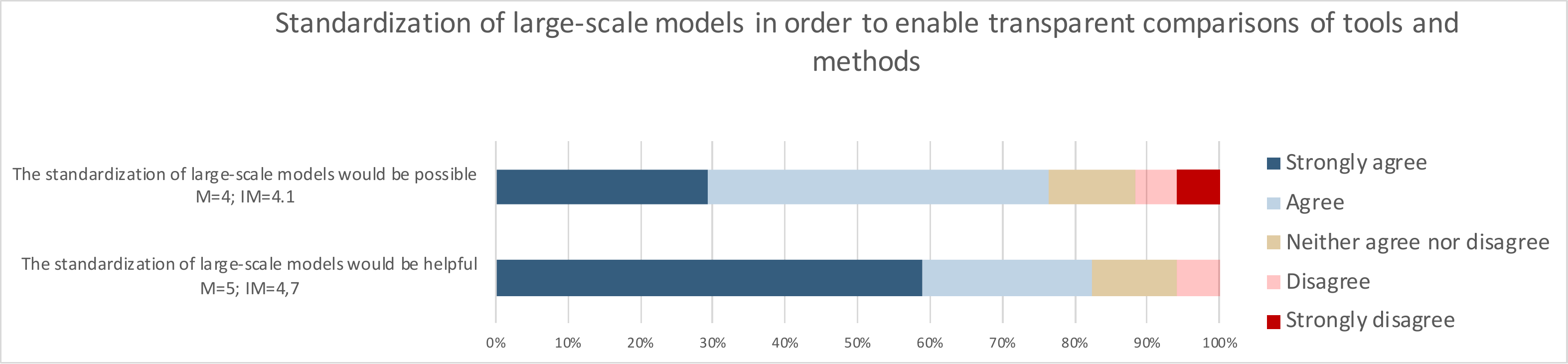}
\caption{Expert assessments: Would it be possible and helpful to define standardized large-scale models in order to enable transparent comparisons of tools and methods 
}
\label{fig:LargeScale}
\end{figure}

\begin{figure}[h!]
\centering
\includegraphics[width=1\textwidth]{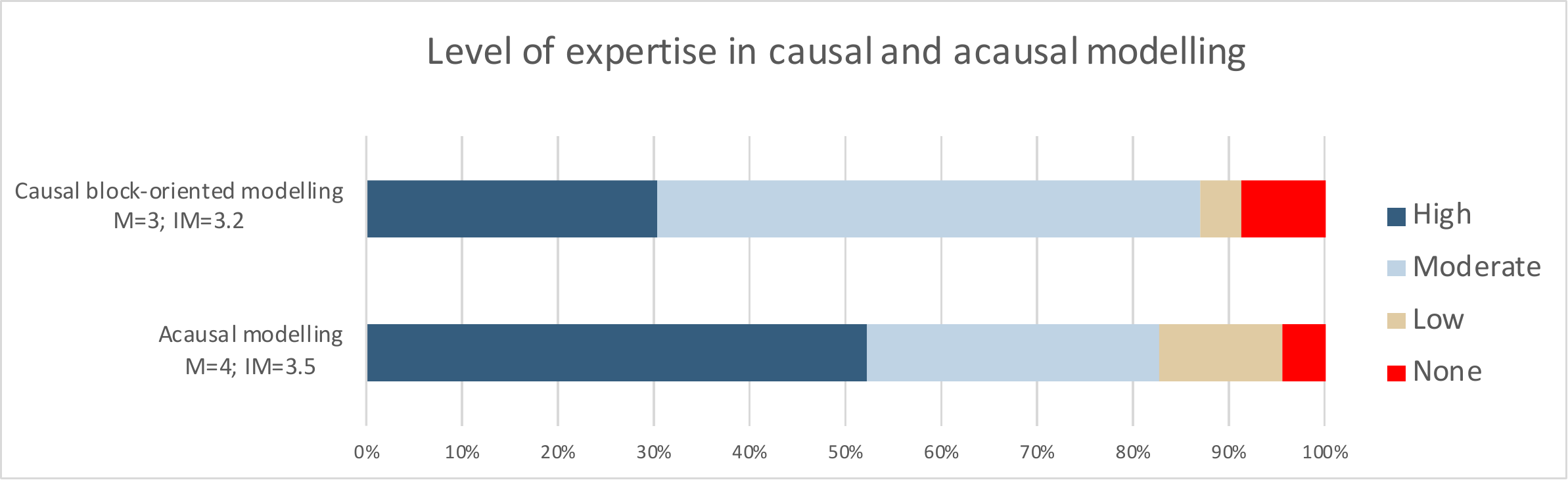}
\caption{Level of expertise in causal and acausal modeling}
 \label{fig:expertise}
 \end{figure}
 
\begin{figure}[h!]
\centering
\includegraphics[width=1\textwidth]{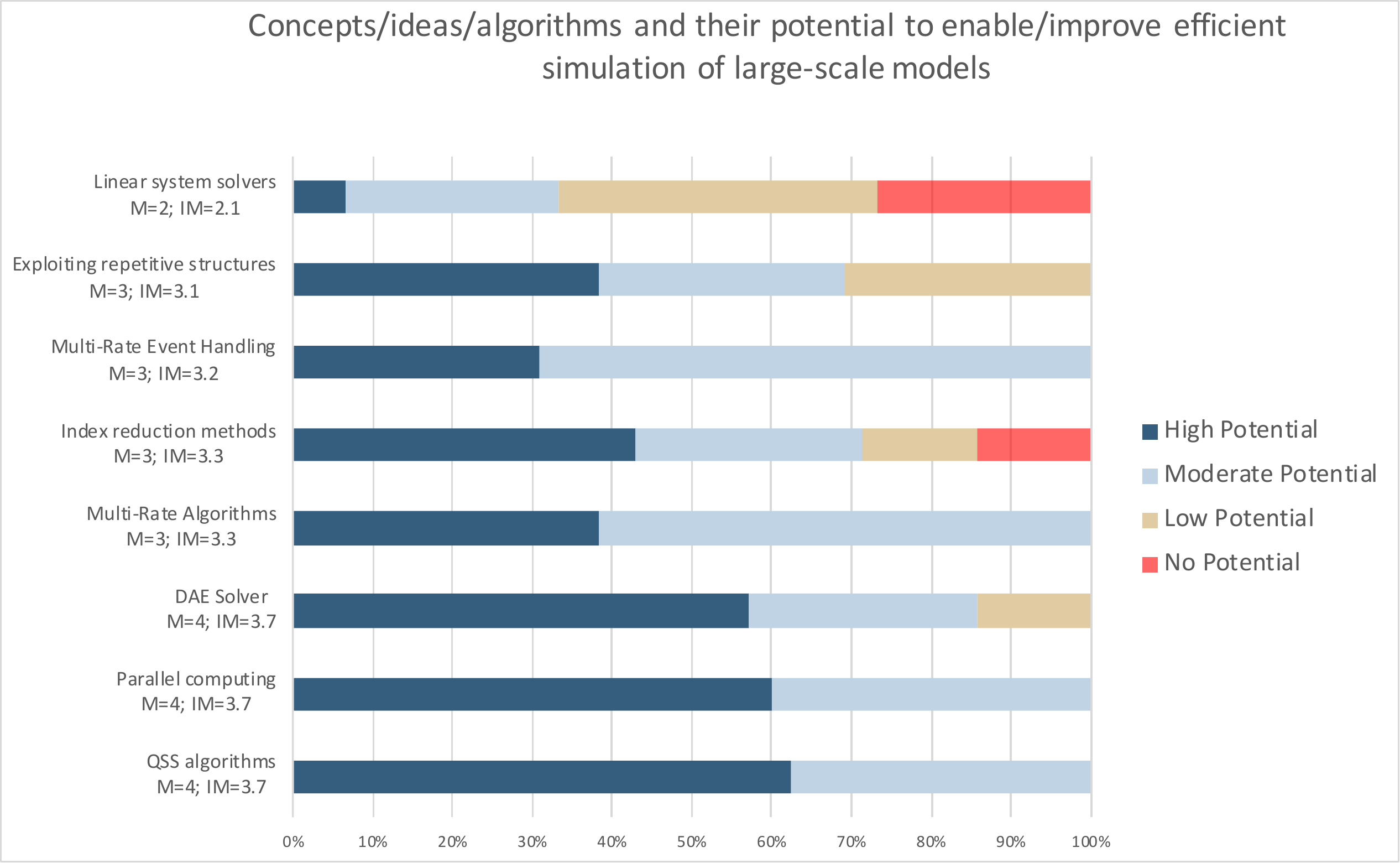}
\caption{Experts were asked to rate concepts/ideas/algorithms and their potential to enable/improve efficient simulation of large-scale models
}
\label{fig:AcausalPotential}
\end{figure}
